%hit.tex: a Plain TeX file by Lucy Martinez and Doron Zeilberger
%How Many Dice Rolls Would It Take to Hit Your Favorite Kind of Number?
%begin macros

\baselineskip=14pt
\parskip=10pt

\magnification=\magstephalf

\def\P{{\cal P}}

\def\1{{\overline{1}}}
\def\2{{\overline{2}}}
\parindent=0pt
\overfullrule=0in

\def\frac#1#2{{#1 \over #2}}
%\headline={\rm  \ifodd\pageno  \RightHead  \else  \LeftHead  \fi}
%\def\RightHead{\centerline{
%Title
%}}
%\def\LeftHead{ \centerline{Doron Zeilberger}}
%end macros
\centerline
{\bf 
How Many Dice Rolls Would It Take to Hit Your Favorite Kind of Number?
}
\bigskip
\centerline
{\it Lucy MARTINEZ and Doron ZEILBERGER}
\bigskip

{\bf Abstract}:  
Noga Alon and Yaakov Malinovsky recently studied the following game: you start at $0$, and keep rolling a fair standard die,  and add the
outcomes until the sum happens to be prime. We generalize this in several ways, illustrating the power of
symbolic, rather than merely numeric, computation. We conclude with polemics why the beautiful rigorous error estimate of Alon and
Malinovsky is only of theoretical interest, explaining why we were content, in our numerous extensions, with non-rigorous, but practically-certain, estimates.

{\bf Preface}

Noga Alon and Yaakov Malinovsky [AM] recently considered the following {\it solitaire} game.

Suppose that you {\it love} prime numbers, and want to reach them by rolling a fair standard die, whose faces are labeled with  $\{1,2,3,4,5,6\}$.
You start at $0$, and you keep adding up the number of dots, and quit as soon as the running total is prime.
They proved, {\it rigorously}, that the expected duration of this `game' is $2.4284\dots$, and that its variance is $6.2427\dots$.
They proceeded in two steps. First they found a candidate approximation (by truncating the implied infinite series), and then used
(a slightly weaker version of) the prime number theorem to rigorously bound the error.

Let's first get a feel for the game by spelling out the first two rounds.

If you are really lucky, the first die-roll is already prime, i.e. you got $2$,$3$ or $5$. So with probability $\frac{3}{6}=\frac{1}{2}$ the `game' only takes
one round. 

Suppose that you did not hit a prime the first time, i.e. you either got $1$, $4$, or $6$.

$\bullet$ If the first throw is $1$, then if the next throw is either $1$, $2$, $4$, or $6$, then you are done, contributing $\frac{1}{6} \cdot \frac{4}{6}=\frac{1}{9}$ .

$\bullet$ If the first throw is $4$, then if the next throw is either $1$, $3$, then you are done, contributing $\frac{1}{6} \cdot \frac{2}{6}=\frac{1}{18}$ .

$\bullet$ If the first throw is $6$, then if the next throw is either $1$, $5$, then you are done, contributing $\frac{1}{6} \cdot \frac{2}{6}=\frac{1}{18}$ .

So the probability that the game lasts exactly two rounds is $\frac{1}{9}+\frac{1}{18}+\frac{1}{18}=\frac{2}{9}$ .

With probability $1-\frac{1}{2}-\frac{2}{9}=\frac{5}{18}$ you need to continue.

If we are really unlucky, with {\bf positive} probability, we may never hit a prime in any prescribed number of moves. 
For example, it is possible, but unlikely, that after getting $1$ followed by $3$, you only get even rolls, so the
partial sums are never odd, let alone prime.

{\bf Question}: How long, on average, would it take until you hit a prime? Or more formally:

{\it What is the {\bf expectation} of the random variable `duration of the prime-seeking' game?}

In the present paper, in addition to many other things, we show that, to  $103$ digits, this real number, that may be named the {\it Alon-Malinovsky} constant, starts with

$2.42849791369350423036165217765842179665120021118534684674615373381696983417849235$

While, unlike the original estimate (with four decimal digits after the decimal point), that is {\bf fully rigorous},
ours is `only' a  {\it non-rigorous estimate}, yet, as we would argue, it is {\bf practically certain}. See the concluding section why this is good enough for us.

{\bf Numerical Dynamical Programming}

The Alon-Malinovsky paper [AM] consisted of two parts. The first was {\it computational}, using {\it dynamical programming}.
They defined the quantity $p(k,n)$, $k \leq n \leq 6k$, where $n$ is a non-prime, to be the
the probability that after $k$ rolls, the running sum is  the non-prime $n$. 

They used the {\bf dynamical programming recurrence}

$$
p(k,n) \, = \, \frac{1}{6} \, \sum_{i} p(k-1,n-i) \quad,
$$
where the sum ranges over all $i$ between $1$ and $6$ so that $n-i$ is a non-prime, to compute many terms, and then defined
$$
p(k+1) := \sum_{\{n: k \leq n \leq 6k\}} p(k,n) \quad,
$$
$$
E_K = \sum_{k=1}^{K} p(k) \quad,
$$
and argued that $E_{1000}$ is a good approximation for the desired expected duration, and then went on to {\it rigorously} bound the error.

We were intrigued, and asked ourselves:

$\bullet$ What if you don't start at $0$ but later on? For example, if you start at $10^{10}$ (obviously a non-prime), how long, on average, would it take until you hit a prime?

$\bullet$ What if instead of a standard die with six faces, you have a different number of faces?

$\bullet$ What if instead of a fair  die  you have a loaded die?

$\bullet$ What if instead of trying to hit a prime you want to hit other kinds of numbers? How long would it take, on average,
to hit a {\it product of two distinct primes}?, {\it product of three distinct primes}?, {\it perfect square} (if you start at a non-square)?, etc. etc.

Our initial approach was to emulate  [AM], namely use dynamical programming, and {\it numerics}, to get very good estimates of the
expected duration, for all these different variants, and we collected lots of data. This is accomplished in the first Maple package accompanying this article, {\tt HIT1.txt}, available from

{\tt https://sites.math.rutgers.edu/\~{}zeilberg/tokhniot/HIT1.txt} \quad .

To get estimates for the expected number of rolls to hit a prime for dice with number of faces from $3$ to $15$, see the output file

{\tt https://sites.math.rutgers.edu/\~{}zeilberg/tokhniot/oHIT1a.txt} \quad .

To see the analogous quantities for reaching {\it product of two distinct primes}, see the output file

{\tt https://sites.math.rutgers.edu/\~{}zeilberg/tokhniot/oHIT1b.txt} \quad .

In particular, it takes, on average, $3.788921291\dots$ rolls to hit a product of two distinct primes with a standard fair die.
Note that this is a {\it non-rigorous} estimate, yet {\bf practically absolutely certain} (see last section).

First we planned also to emulate the rigorous error-analysis in [AM] for these other scenarios, but got an {\it epiphany},
it is not worth the trouble! See the last section why.

Then we got a second epiphany. Using {\it symbolic computation} (with Maple) rather than {\it numeric computation} (with Matlab) to handle this kind of problems is more natural and streamlined, and possibly more efficient.

{\bf Using Symbolic Computation to Model the Game}

Sooner or later (with probability $1$) the game ends, at some number of rounds, when you reached a certain prime.
Let $q(k,n)$ be the probability that it ended after $k$ rounds and that the running sum then was the prime $n$.

Everything about this process is encoded in the bivariate {\it probability generating function}, the {\it infinite} double-series
$$
F(t,x) \, := \, \sum_{k=1}^{\infty} \left ( \sum_{ {{k \leq n <6k} \atop { n \,\, prime}}} q(k,n) x^n \right ) t^k \quad .
$$

Of course $F(1,1)=1$, and the expected duration is $F_t(1,1)$, while the expected final location is $F_x(1,1)$. The variance, higher moments,
and {\it mixed moments} (in particular the {\it covariance} [from which we get the more informative {\it correlation}]) could be gotten by differentiating with respect to $t$ and/or $x$ and then substituting $x=1,t=1$.

Alas, this is an {\it infinite series}, so let's be more modest and try and compute the {\bf truncated series}, for a given {\bf finite} maximal number of rounds, $R$:
$$
F_R (t,x) \, := \, \sum_{k=1}^{R} \left ( \sum_{ {{k \leq n <6k} \atop { n \,\, prime}}} q(k,n) x^n \right ) t^k \quad .
$$

To illustrate how we got Maple  to compute $F_R(t,x)$, we will continue with the original game of starting at $0$, rolling a fair standard die, and seeking a prime.
The same approach works in general, and that is what we implemented.

Let $P(x)$ be the {\bf probability generating function} of the die:
$$
P(x)=\frac{1}{6} \sum_{i=1}^{6} x^i \,= \, \frac{1}{6}\,x \,+ \,\frac{1}{6}\,x^2 \,+ \,\frac{1}{6}\,x^3 \,+ \,\frac{1}{6}\, x^4+ + \,\frac{1}{6}\, x^5 \, + \, \,\frac{1}{6}\, x^6 \quad .
$$

We need the following {\it operator} defined on polynomials $\sum_{i=1}^{n} a_i x^i$
$$
\P \left ( \sum_{i=1}^{n} a_i x^i \right ) \, := \,\sum_{{{ 1 \leq i \leq n} \atop {i \,\, prime}}} a_i x^i  \quad.
$$

For example
$$
\P(x+3x^2+5x^4+\frac{1}{2} x^6+ 8x^7) = 3x^2+ 8x^7 \quad.
$$

We also need an auxiliary sequence of polynomials $S_R(x)$, that takes care of the {\it survivors} at the $R^{th}$ round.

Initialize: $S_0(x):=1$ (more generally $S_0(x):=x^{init}$).  Also $F_0(t,x):=0$.

Suppose that you already have $F_{R-1}(t,x)$.

If currently you are at the $R$-th round, with the   previous survival polynomial, $S_{R-1}(x)$, define
$$
N_R(x):=\P(P(x) S_{R-1}(x)) \quad, \quad S_R(x):=P(x)S_{R-1}(x)-N_R(x) \quad F_R(t,x) := F_{R-1}(t,x)+ N_R(x)\,t^R \quad .
$$

Let's explain. $P(x)\,S_{R-1}(x)$ is the probability generating function, according to location, of the {\it new guys}, some of them prime, and some not. Applying $\P$ extracts
the primes, the  new {\it inductees}.

Let's illustrate the first two steps. At the first round:
$$
N_1(x)=\P(1 \cdot  (\frac{1}{6}\,x \,+ \,\frac{1}{6}\,x^2 \,+ \,\frac{1}{6}\,x^3 \,+ \,\frac{1}{6}\, x^4 \,+  \,\frac{1}{6}\, x^5 \, + \,\frac{1}{6}\, x^6))=
\,\frac{1}{6}\,x^2 \,+ \,\frac{1}{6}\,x^3 \,+ \,\frac{1}{6}\, x^5 \quad .
$$
So 
$$
S_1(x)= \,\frac{1}{6}\,x \,+ \,\frac{1}{6}\,x^4 \,+ \,\frac{1}{6}\, x^6 \quad .
$$
and
$$
F_1(t,x)=\left (\frac{1}{6} x^{2}+\frac{1}{6} x^{3}+\frac{1}{6} x^{5} \right ) t \quad .
$$
 Next:
$$
 S_1(x) \, P(x) = \, (\frac{1}{6}\,x \,+ \,\frac{1}{6}\,x^4 \,+ \,\frac{1}{6}\, x^6)
\, (\frac{1}{6}\,x \,+ \,\frac{1}{6}\,x^2 \,+ \,\frac{1}{6}\,x^3 \,+ \,\frac{1}{6}\, x^4+ + \,\frac{1}{6}\, x^5 \, + \, \,\frac{1}{6}\, x^6)
$$
$$
= \, \frac{1}{36} x^{2}+\frac{1}{36} x^{3}+\frac{1}{36} x^{4}+\frac{1}{18} x^{5}+\frac{1}{18} x^{6}+\frac{1}{12} x^{7}+\frac{1}{18} x^{8}+\frac{1}{18} x^{9}+\frac{1}{18} x^{10}+\frac{1}{36} x^{11}+\frac{1}{36} x^{12} \quad .
$$
Applying $\P$ we get
$$
N_2(x)=\frac{1}{36} x^{2}+\frac{1}{36} x^{3}+\frac{1}{18} x^{5}+\frac{1}{12} x^{7}+\frac{1}{36} x^{11} \quad.
$$
So,
$$
F_2(t,x)=\left(\frac{1}{6} x^{2}+\frac{1}{6} x^{3}+\frac{1}{6} x^{5}\right) t +
\left(\frac{1}{36} x^{2}+\frac{1}{36} x^{3}+\frac{1}{18} x^{5}+\frac{1}{12} x^{7}+\frac{1}{36} x^{11}\right) t^{2} \quad ,
$$
$$
S_2(x)=\frac{1}{36} x^{4}+\frac{1}{18} x^{6}+\frac{1}{18} x^{8}+\frac{1}{18} x^{9}+\frac{1}{18} x^{10}+\frac{1}{36} x^{12} \quad,
$$
and we keep going until we reach the $R^{th}$ round.

The probability that the game ends in $\leq R$ rounds is $F_R(1,1)$, (for $R$ large this is very close to $1$). Also of interest is
the {\bf conditional} probability generating function
$$
{\overline F}_R(t,x) :=\frac{F_R(t,x)}{F_R(1,1)} \quad .
$$

{\bf The Maple package {\tt HIT2.txt}}

This symbolic-computational approach is implemented in the Maple package {\tt HIT2.txt} available from:

{\tt https://sites.math.rutgers.edu/\~{}zeilberg/tokhniot/HIT2.txt} \quad .

To get a list of the main procedures type {\tt ezra();}. For example procedure

{\tt GFpG(n,R,init,t,x,P)}

inputs:

$\bullet$ a positive integer {\tt n}, corresponding to the number of faces in a fair die ;

$\bullet$ a positive integer {\tt R}, corresponding to the maximum number of rolls ;

$\bullet$ a non-negative integer  {\tt init}, the starting location ;

$\bullet$  (formal) variables {\tt t,x} ;

$\bullet$  a property {\tt P} (e.g. {\tt isprime}) .

It outputs the truncated bivariate probability generating function up to $t^R$,
(what we called above $F_R(t,x)$).  For example to get $F_{100}(t,x)$ for the original [AM] scenario, type:

{\tt GFpG(6,100,0,t,x,isprime);} \quad .

{\bf Data}

Once we have taken the trouble to write the Maple code, we can generate lots of interesting data.

$\bullet$ If you want to see the number of rounds it takes to guarantee that you reached a prime, starting at $0$,  with probability $\geq 1-10^{-7}$, as well as
the {\it expected duration}, {\it variance}, {\it skewness}, and {\it kurtosis} for dice with number of faces from $2$ to $40$ look here:

{\tt https://sites.math.rutgers.edu/\~{}zeilberg/tokhniot/oHIT2a.txt} \quad .

$\bullet$ If you want to see the number of rounds it takes to guarantee that you reached a prime, starting at $0$,  with probability $\geq 1-10^{-20}$, as well as
the {\it expected duration}, {\it variance}, {\it skewness}, and {\it kurtosis} for dice with number of faces from $2$ to $40$ look here:

{\tt https://sites.math.rutgers.edu/\~{}zeilberg/tokhniot/oHIT2a1.txt} \quad .

$\bullet$ If you want to see the number of rounds it takes to guarantee that you reached a {\bf product of two distinct primes}, starting at $0$,  with probability $\geq 1-10^{-7}$, as well as
the {\it expected duration}, {\it variance}, {\it skewness}, and {\it kurtosis} for dice with number of faces from $2$ to $40$ look here:

{\tt https://sites.math.rutgers.edu/\~{}zeilberg/tokhniot/oHIT2b.txt} \quad .

Note in particular that for a standard (six-faced) fair die,  the expected duration is $3.7889\dots$, a bit longer than for hitting a prime.

$\bullet$ If you want to see the number of rounds it takes to guarantee that you reached a {\bf product of three distinct primes}, starting at $0$,  with probability $\geq 1-10^{-7}$, as well as
the {\it expected duration}, {\it variance}, {\it skewness}, and {\it kurtosis} for dice with number of faces from $2$ to $40$ look here:

{\tt https://sites.math.rutgers.edu/\~{}zeilberg/tokhniot/oHIT2c.txt} \quad .

Note in particular that for a standard (six-faced) fair die,  the expected duration is  $17.616887\dots$, quite a bit longer than for reaching a product of two distinct primes.

$\bullet$ If you want to see the number of rounds it takes to guarantee that you reached a {\bf product of four distinct primes}, starting at $0$,  with probability $\geq 1-10^{-7}$, as well as
the {\it expected duration}, {\it variance}, {\it skewness}, and {\it kurtosis} for dice with number of faces from $2$ to $40$ look here:

{\tt https://sites.math.rutgers.edu/\~{}zeilberg/tokhniot/oHIT2d.txt} \quad .

Note in particular that for a standard (six-faced) fair die, the expected duration is  $112.907872\dots$ much longer than for reaching a product of three distinct primes.

$\bullet$ If you want to see the number of rounds it takes to guarantee that you reached a {\bf perfect square}, starting at $2$ (of course $0$ and $1$ are perfect squares),  with probability $\geq 1-10^{-6}$, as well as
the {\it expected duration}, {\it variance}, {\it skewness}, and {\it kurtosis} for dice with number of faces from $2$ to $40$ look here:

{\tt https://sites.math.rutgers.edu/\~{}zeilberg/tokhniot/oHIT2e.txt} \quad .

Note in particular that for a standard (six-faced) fair die,  the expected duration is $9.01861\dots$.

$\bullet$ Suppose that you allow up to 200 dice-rolls (don't worry, the probability that you won't get a prime by the  time you roll at most $200$ times is less than $10^{-18}$ in all 
the cases here),
 to see not only the estimated expected duration, but also the expected destination, as well as the correlation 
(not surprisingly close to $1$, but a bit surprisingly not that close, e.g. for the usual six-faced fair die it is $0.965644\dots$).
For fair dice with number of faces from $2$ to $20$  look here:

{\tt https://sites.math.rutgers.edu/\~{}zeilberg/tokhniot/oHIT2g.txt} \quad .

$\bullet$ So far we only treated  {\it fair} dice, but our Maple code can equally well treat {\it loaded} dice. For one example, see

{\tt https://sites.math.rutgers.edu/\~{}zeilberg/tokhniot/oHIT2k.txt} \quad .

$\bullet$ To see nice plots how, with various starting places, the expected duration changes with the number of faces, see

{\tt https://sites.math.rutgers.edu/\~{}zeilberg/tokhniot/oHIT2pics.pdf} \quad .

{\bf How to Estimate (Non-Rigorously but Practically Certainly) the Alon-Malinovsky Constant and Why it is not So Interesting}

Let's face it. Life is finite. Also, as much as you love primes, it would be very tedious to keep rolling a die. So beforehand you decide
what is the maximum number of rolls that you are willing to make. Once you decide {\it beforehand} about the number of maximum 
rolls, $R$, you would like to know:

$\bullet$ What is the probability that you would indeed achieve your goal of getting a prime in $\leq R$ rolls? Let's call it $a_R$.

$\bullet$ {\it Conditioned} on that event, what is the expected number of rolls that it will take to finish? Let's call it $M_R$.

$\bullet$ {\it Conditioned} on that event, what is the expected location (relative to the starting place) where you wind up at? Let's call it $L_R$.

Of course you can also ask about the variance, and higher moments, and even mixed moments.

To compute these quantities, you first find the bi-variate truncated probability generating function, $F_R(t,x)$, that our Maple package {\tt HIT2.txt} computes using

{\tt GFpG(n,R,init,t,x,P);} (see above).

For the original [AM] case we have  {\tt n=6}, {\tt init=0} and  {\tt P=isprime}. So the function call to get our $F_R(t,x)$ is:

{\tt GFpG(6,R,0,t,x,isprime);}   \quad,

for any desired positive integer {\tt R}.

We have
$$
a_R=F_R(1,1) \quad.
$$
Define, as above
$$
{\overline F}_R(t,x) :=\frac{F_R(t,x)}{a_R} \quad,
$$
the conditional probability generating function, conditioned on terminating in $\leq R$ rolls.

We have
$$
M_R=\frac{\partial}{\partial t} {\overline F}_R(t,x) \vert_{x=1,t=1} \quad,
$$
for the conditional expected duration (conditioned on finishing in $\leq R$ rolls), and
$$
L_R=\frac{\partial}{\partial x} {\overline F}_R(t,x) \vert_{x=1,t=1} \quad ,
$$
for the conditional expected exiting location.

Using our Maple program with $R=200$, we get
$$
a_{200}=1-2.9020152044089\cdot 10^{-19} \quad.
$$
In other words, the probability that you would have to roll more than $200$ rolls is minuscule. Assuming that you indeed finished in $\leq 200$ rolls, we have
$$
M_{200}=2.4284979136935041712\dots \quad,
$$
$$
L_{200}=8.49974269792726459237146481486\dots \quad .
$$

Using our Maple program with $R=400$, we get
$$
a_{400}=1- 1.32546541967224185265621962\cdot 10^{-33} \quad .
$$
In other words, the probability that you would have to roll more than $400$ rolls is even more minuscule. Assuming that you indeed finished in $\leq 400$ rolls, we have
$$
M_{400}=2.4284979136935042303660819062417645\dots \quad ,
$$
$$
L_{400}=8.4997426979272648062812866718461364\dots  \quad .
$$

Using our Maple program with $R=1000$, we get
$$
a_{1000}=1-2.183194254589149\cdot 10^{-73} \quad .
$$
In other words, the probability that you would have to roll more than $1000$ rolls is much less than all of us dying in a nuclear holocaust!
Assuming that you indeed finished in $\leq 1000$ rolls we have
$$
M_{1000}=2.428497913693504230366081906241764513835\dots \quad,
$$
$$
L_{1000}=8.4997426979272648062812866718480475 \quad.
$$

Note that these values are {\bf exact}, and Maple has them as {\bf rational numbers}, and by resetting {\tt Digits}, one can get as many decimals as one wishes.

Now the Alon-Malinovsky constant may be defined by
$$
M_{\infty} := \lim_{R->\infty} M_R \quad,
$$
and analogously
$$
L_{\infty} := \lim_{R->\infty} L_R \quad.
$$

Since $|M_{400}-M_{1000}|\leq 10^{-31}$, i.e. they agree to the first 30 digits, but to {\bf play it safe}, let's only take the first $20$ digits of $M_{1000}$ (or for that matter,
the first $20$ digits of $M_{400}$, since they are identical) and we have the following:

{\bf A Non-rigorous but Practically Certain Estimate} of the Alon-Malinovsky constant is: $2.4284979136935042304\dots$. 

Taking larger $R$ we get the more precise estimate, to $103$ digits, stated at the beginning of this article. See the output file

{\tt https://sites.math.rutgers.edu/\~{}zeilberg/tokhniot/oHIT2h.txt } \quad .

But, {\it frankly}, rigorous or not,  $M_{\infty}$  is {\bf not that interesting!} Life is finite. The {\it assurance} that with probability larger than $1-10^{-72}$ you will finish the game in $\leq 1000$ rolls,
and conditioned on that, the expected number of rolls, i.e. $M_{1000}$, is $2.428497913693504230366081906241764513835\dots$ is much more interesting, and {\it useful}!
Let's face it, life is finite, and even during our short life, we have better things to do than roll a die until we get a prime. The good news is that with
very high probability we will get there soon enough!

{\bf References}

[AM] Noga Alon and Yaakov Malinovsky, {\it Hitting a prime in 2.43 dice rolls (on average)},    arXiv:2209.07698 [math.PR], 16 Sep 2022 (v1), 21 Dec 2022 (v3). \hfill\break
{\tt https://arxiv.org/abs/2209.07698} \quad .

\bigskip
\hrule
\bigskip
Lucy Martinez and Doron Zeilberger, Department of Mathematics, Rutgers University (New Brunswick), Hill Center-Busch Campus, 110 Frelinghuysen
Rd., Piscataway, NJ 08854-8019, USA. \hfill\break
Email: {\tt  lm1154 at math dot rutgers dot edu} \quad, \quad {\tt DoronZeil at gmail dot com}   \quad .

Written: {\bf Jan.  31, 2023}.

\end